\documentclass{amsart} 
\usepackage{amssymb, amscd}
\numberwithin{equation}{section}

\def\CC{{\mathbb C}}  

\def\EE{{\mathbb E}}    
\def\FF{{\mathbb F}}  
\def\GG{{\mathbb G}}

\def\PP{{\mathbb P}}
\def\QQ{{\mathbb Q}}

\def\ZZ{{\mathbb Z}}

\def\Acal{{\mathcal A}}

\def\Fcal{{\mathcal F}}

\def\Mcal{{\mathcal M}}

\def\Ocal{{\mathcal O}}
\def\Pcal{{\mathcal P}}

\def\Xcal{{\mathcal X}}

\newcommand\Tr{\operatorname{Tr}}
\renewcommand\Sp{\operatorname{Sp}}

\let\co\colon
\newcommand\Definition[1]{\emph{#1}}
\newcommand\map[3]{\ensuremath{{#1}\co{#2}\to{#3}}}

\newcommand\proofsquare{\nobreak\hfill \hbox{%
\vrule height 5pt 
\kern-.4pt
 \vbox{%
\hrule width 5pt depth0pt height.4pt
 \kern4.6pt \hrule  }
\kern-3.75pt 
\vrule height 5pt}\kern1pt
\par}

\newtheorem{theorem}{Theorem}[section]
\newtheorem{lemma}[theorem]{Lemma}
\newtheorem{proposition}[theorem]{Proposition}
\newtheorem{corollary}[theorem]{Corollary}

\newtheorem{definition-lemma}[theorem]{Definition-Lemma}

\theoremstyle{definition}

\newtheorem{example}[theorem]{Example}

\theoremstyle{remark} 
\newtheorem{remark}[theorem]{Remark}

\begin{document}

\title[Top Chern class of the Hodge bundle]{The Top Chern Class of the Hodge 
Bundle\\
on the Moduli Space of Abelian Varieties
\footnote{\tt topchernclass.tex\today}} 
\author{Torsten Ekedahl}
\address{Department of Mathematics\\
 Stockholm University\\
 S-106 91  Stockholm\\
Sweden}
\email{teke@matematik.su.se}
\author{Gerard van der Geer}

\address{Faculteit Wiskunde en Informatica, University of
Amsterdam, Plantage Muidergracht 24, 1018 TV Amsterdam, The Netherlands}

\email{geer@wins.uva.nl} 


\begin{abstract}

We give upper and lower bounds for the order of the top Chern class 
of the Hodge bundle on the moduli space of principally
polarized abelian varieties.
We also give a generalization to higher genera 
of the famous formula $12 \, \lambda_1=\delta$ for genus $1$. \end{abstract}

\maketitle

\begin{section}{Introduction}
\label{sec: intro}
\bigskip
\noindent
Let ${\Acal}_g/ \ZZ$ denote the moduli stack of principally polarized
abelian varieties of dimension~$g$. This is an irreducible algebraic
 stack of relative dimension $g(g+1)/2$ with irreducible fibres over
$\ZZ$. The stack ${\Acal}_g$  carries a locally free sheaf
$\EE$ of rank $g$, the Hodge bundle, defined as follows.
If $A/S$ is a principally polarized
abelian variety with zero section $s$  we get a locally free
sheaf $s^*\Omega^1_{A/S}$
of rank $g$ on $S$ and this is compatible with pull backs. If $\pi: A \to
S$ denotes  the structure map it satisfies the property
$\Omega^1_{ A/S}= \pi^*(\EE)$. The Hodge bundle can be
extended to a locally free sheaf (again denoted by) $\EE$ on every
smooth toroidal compactification $\tilde {\Acal}_g$ of ${\Acal}_g$,
cf.\  \cite{F-C}.

The Chern classes $\lambda_i$ of the Hodge bundle $\EE$ are defined over $\ZZ$
and give for each fibre ${\Acal}_g\otimes k$ rise to classes $\lambda_i$ in the
Chow ring $CH^*({\Acal}_g\otimes k )$, and in $CH^*(\tilde{\Acal}_g\otimes
k)$. They generate subrings ($\QQ$-subalgebras) of $CH_{\QQ}^*({\Acal}_g\otimes
k)$ and of $CH_{\QQ}^*(\tilde {\Acal}_g\otimes k)$ which are called the
\Definition{tautological subrings}.

 It was proved in \cite{vdG1} by an application of the
Grothendieck-Riemann-Roch theorem that these classes in the Chow ring
$CH_{\QQ}^*( {\Acal}_g)$ with rational coefficients satisfy the
following relation 
\begin{equation}\label{fund rel}
(1+\lambda_1+\ldots +\lambda_g)(1-\lambda_1+\ldots +(-1)^g
\lambda_g)=1. %\eqno(1)
\end{equation}
Furthermore, it was proved that $\lambda_g$ vanishes in the Chow group
$CH_{\QQ}({\Acal}_g)$ with rational coefficients.  
The class $\lambda_g$ does not vanish on $\tilde{\Acal}_g$.
This raises two questions.  First, since $\lambda_g$ is a torsion class
on ${\Acal}_g$ we may ask what its order is. Second, since
$\lambda_g$ up to torsion comes from a class on the `boundary'
$\tilde{\Acal}_g - {\Acal}_g$   we may ask for a description of this
class. As an answer to these questions we give an upper  bound on the order of
$\lambda_g$ in the second section and a non-vanishing result
in the third section which implies a lower bound.
 We then give a several  representing cycles for
the top Chern class: on ${\Acal}_g$,  on the Satake compactification,
and  in the last
section we describe a cycle class in a partial
compactification of
${\Acal}_g$ which represents $\lambda_g$ in the Chow group $CH^g_{\QQ}$.
This result can be viewed  as a generalization 
of the well-known relation $12\, \lambda_1=\delta$
for $g=1$.\bigskip
\end{section}

\begin{section}{A bound on the order of the class $\lambda_g$.}

We will make some computations in the Chow group of $\Acal_g$. For the
definition of the (integral) Chow groups for stacks we refer to \cite{Kr}.
We begin with a lemma which is no doubt well known.
\begin{lemma}\label{chern}
If $E$ is a vector bundle of rank $g$, the total Chern class of the
graded vector bundle $\Lambda^*E$ is zero in degrees 1 to $g-1$ and
$-(g-1)!c_g$ in degree $g$.
\end{lemma}
%% \noindent
%% {\sl Proof.} 
\begin{proof}
As usual we note that the components of the total Chern
class are universal polynomials in the Chern classes
$c_i:=c_i(E)$ and we may let $E$
be the universal bundle. Then as we may think of the coefficients in
the universal polynomials as rational numbers we can note that the
Chern classes of degree 1 to $g-1$ of
 $\Lambda^*E$ vanish if and only if the Newton polynomials of the same
degrees do, that is, if and only if ${\text ch}(\Lambda^*E)$ vanishes in the
same degrees. Thus the first part followed from the Borel-Serre formula
(cf.\ \cite{B-S})
$$
{\text ch}(\Lambda^*E)=(-1)^gc_g{\text{Td}}(E)^{-1}.
$$
Furthermore, if the Chern classes of degree 1 to $g-1$ vanish for a
bundle $F$, then it is clear that $(-1)^{g-1}g\, c_g(F)=s_g(F)$, as can be
seen from Newton's formula. Hence in degree $g$ we have ${\text ch}(\Lambda^*E)$ is
$(-1)^{g-1}gc_g(F)/g!$ and using the Borel-Serre formula again gives
the desired formula.
\end{proof}
\begin{lemma}\label{orderp}
Let $p$ be a prime. Suppose that $\pi: A \to S$ is a family of abelian varieties of
relative dimension $g$, where $S$ is a scheme, and $L$ is a line bundle on $A$
of order $p$, on all fibres of $\pi$. If $p>\min(2g,\dim S+g)$ then
$p(g-1)!\lambda_g=0$.
\end{lemma}
%% \noindent
%% {\sl Proof.}  
\begin{proof}
By twisting $L$ by a line bundle on $S$ so that it is trivial along the zero
section we may assume that it is of order $p$ on $A$.
Denoting the class of $L$ in $K_0$ by $[L]$ we then either have
that $p([L]-1)$ has support of codimension $>2g$ if $p>2g$ or is zero if $p>\dim
S+g$. Indeed, from the relation
$$0=[L]^p-1=p([L]-1)(1+(p-1)/2([L]-1)+...)+([L]-1)^p$$ and the fact that
 $[L]-1$ is nilpotent we get that $p([L]-1)$ is divisible by
$([L]-1)^p$ and that element is supported in codimension $\ge p$. Now in
the first case the image under
${\pi}$  of the support has codimension $>g$ on $S$ and so we may safely
remove it and may assume that $p[L]=p$ in $K_0(A)$.

Consider now the Poincar\'e bundle ${\Pcal}$ on $A\times_S\check A$, where
$\check A$ denotes the dual abelian variety.  By base change $R\pi_*{\Ocal}_A$
is the (derived) pullback along the zero-section of $\check A$ of the sheaf
$R\pi_*\Pcal$. We have that $p[{\Pcal}]=p[L\bigotimes {\Pcal}]$ and so
$p[R\pi_*{\Pcal}]=p[R\pi_*(L\bigotimes {\Pcal})]$. Now, a fibrewise calculation
shows that $R\pi_*(L\bigotimes \Pcal)$ has support along the inverse of the
section of $\check A$ corresponding to $L$. As that section is everywhere
disjoint from the zero section the pullback of $R\pi_*(L\bigotimes \Pcal)$ along
the zero section is 0 and thus $p[R\pi_*{\Ocal}_A]=0$. Using lemma \ref{chern}
and applying the total Chern class then gives $1=(1-(g-1)!\lambda_g+\ldots)^p$
which in turn gives the lemma.
\end{proof}
\begin{definition-lemma}\label{lcd}
For an integer $g$ we let  $n_g$ be  the largest common divisor of
all $p^{2g}-1$ where $p$ runs through all primes larger than a sufficiently large fixed number $N$
(which may be taken to be $2g+1)$). For an odd prime $p$, the exact power $p^k$
of $p$ that divides $n_g$ is the largest $k$ such that $p^{k-1}(p-1)$ divides
$2g$ and $0$ if $p-1$ does not divide $2g$. The exact power $2^k$ that divides
$n_g$ is the largest $k$ such that $2^{k-2}$ divides $2g$.
\end{definition-lemma}
\begin{proof}
This follows directly from the structure of $(\ZZ/p^k)^*$ and Dirichlet's prime
number theorem.
\end{proof}
\smallskip
\begin{example}
We have $n_1=24$, $n_2=240$, $n_3=504$ and $n_4=480$.
\end{example}

\smallskip
\begin{lemma} We have 
$$ 
\prod_{i=1}^g n_i = \prod_{p}
([\frac{2gp}{p-1}]!)_p,
$$
where $p$ runs over the primes. 
\end{lemma}
\smallskip

Define for $g\geq 1$  the positive rational number 
$$
 p(g):= (-1)^g\prod_{j=1}^g \frac{\zeta(1-2j)}{2},
$$ 
where $\zeta(s)$ is the Riemann zeta-function. By the Proportionality
theorem of Hirze\-bruch-Mumford and a theorem of Siegel-Harder we know that
the degree of $\prod_{j=1}^g
\lambda_j$ equals
$\# {\rm Sp}(2g,\ZZ/n)\times p(g)$ on the scheme ${\Acal}_g[n]$ of abelian varieties
with level $n$ structure.  For
$n\geq 3$ this degree must be an integer. This implies the following corollary.

\smallskip
\begin{corollary} The rational number $\prod_{i=1}^g n_i$  is a multiple of
 the denominator of $p(g)$.
\end{corollary}

\begin{proposition}\label{lambdascheme}
Suppose that $\pi: A\to S$ is a family of abelian varieties of
relative dimension $g$, where $S$ is a scheme. Then $(g-1)!n_g\lambda_g=0$
on $S$.
\end{proposition}
%% \noindent
%% {\sl Proof.} 
\begin{proof}
For any prime $p$ larger than $2g$ we can apply  lemma \ref{orderp} on the
cover obtained by adding a line bundle everywhere of order $p$. Projecting down
to $S$ again and using that that cover has degree $p^{2g}-1$ gives
$(g-1)!p(p^{2g}-1)\lambda_g=0$. We then finish by using Definition \ref{lcd} (and
noting that the factor $p$ causes no trouble as by using several primes we
see that no prime $>2g$ can divide the smallest annihilating integer).
\end{proof}

\begin{proposition}\label{lambdastack}
Suppose that $\pi: A\to S$ is the universal family of abelian varieties of
relative dimension $g$. Then $(g-1)!\prod_{i=1}^gn_i\lambda_g=0$ on $S$.
\end{proposition}
%% \noindent
%% {\sl Proof.} 
\begin{proof}
The proof is almost the same as that of \ref{lambdascheme} only
that now we consider instead the cover given by putting a full level
$p$-structure on $A$. This time we therefore get that
$(g-1)!p^m(p^{2g}-1)(p^{2(g-1)}-1)\ldots\lambda_g$, where $m$ is some irrelevant
positive integer. Using again \ref{lcd} we conclude.
\end{proof}
\begin{example}
i) For $g=1$ we get $24\lambda_1=0$ which is off by a factor 2.

ii) For $g=2$ we get $24\cdot(16\cdot3\cdot5)\lambda_2=0$.

iii) For $g=3$ we get
$2\cdot24\cdot(16\cdot3\cdot5)\cdot(8\cdot9\cdot7)\lambda_3=0$.
\end{example}
\end{section}

%;;%%de Rham
\begin{section}{The order of the Chern Classes of the de Rham Bundle}

We will now consider the Chern classes of the bundle of first relative de Rham
cohomology of the universal abelian variety over $\Acal_g$. We will apply our
results to obtain information on the order of $\lambda_g$ but the results we
obtain should be of independent interest. Note that the \Definition{de Rham
bundle}, $H_{dR}^1:= R^1\pi_*\Omega^\cdot_{{\Xcal}_g/{\Acal}_g}$, where
\map\pi{{\Xcal}_g}{{\Acal}_g} is the universal family, is provided with an
integrable connection and hence its Chern classes in integral or $\ell$-adic
cohomology\footnote{$\ell$ of course being a prime different from the
characteristic}, which we will denote $r_i$, are torsion classes. In this
section we will determine their exact order.

We begin by using a result of Grothendieck to get an upper bound for the order
of $r_i$.
%;;%%%Upper bound for order Proposition upper bound
\begin{proposition}\label{upper bound}
i) We have $r_i=0$ for odd $i$.

ii) We have $n_ir_{2i}=0$.
\end{proposition}
%% {\sl Proof.}
\begin{proof}
The first part follows immediately because $H_{dR}^1$ is a symplectic vector
bundle. 

As for the second part we may assume that the characteristic is $0$ as the case
of positive characteristic follows from the characteristic $0$ case by
specialisation. In that case we may further reduce to the case of the base field
being the complex number. We may also prove the annihilation of $n_ir_{2i}$ in
$\ell$-adic cohomology for a specific (but arbitrary) prime $\ell$.

Now, the existence of Gauss-Manin connection on $H^1_{dR}$ means that it has a
discrete structure group. More precisely, the fundamental group of the algebraic
stack ${\Acal}_g$ is $\Sp_{2g}(\ZZ)$ and $H^1_{dR}$ is the vector bundle
associated to the representation of it given by the natural inclusion of
$Sp_{2g}(\ZZ)$ in $Sp_{2g}(\CC)$. This complex representation is obviously
defined over the rational numbers  so we may apply \cite[4.8]{Gr} with field of
definition $\QQ$. We thus conclude that $r_{2i} \in
H^{2i}({\Acal}_g,\ZZ_\ell(i))$ is killed $\ell^{\alpha(i)}$, where $\alpha(i)$
is defined as
\begin{displaymath}
\inf_{\lambda \in H}v_\ell(\lambda^i-1)
\end{displaymath}
and $H \subseteq \ZZ_{\ell}^\ast$ is the image of the Galois group of the field
of definition of the cyclotomic character. However, as the base field is $\QQ$
this image is all $\ZZ_{\ell}^\ast$ and the result follows from the definition
of $n_i$.
\end{proof}
We now aim to show that this upper bound is the precise order of the $r_i$. We
will do this by pulling them back to classifying spaces for certain finite
abelian groups in whose cohomology we will be able to determine their
images. Over the complex numbers this can be done directly by mapping these
finite groups into $Sp_{2g}(\ZZ)$ and using that the cohomology of $\Acal_g$
equals the cohomology of $Sp_{2g}(\ZZ)$. In the positive characteristic case, if
one knew that the specialisation map in the cohomology of $\Acal_g$ induced an
isomorphism then the positive characteristic result would reduce to the
characteristic $0$ one. Though such a specialisation result seems rather
straightforward using the proper and smooth base change theorems, the toroidal
compactifications of Chai and Faltings and an induction on $g$ we know of no
reference and rather than carrying such an argument through will use another
argument. As a motivation for that argument let us begin by noting that any
finite subgroup $G$ of $Sp_{2g}(\ZZ)$ arises as the automorphism group of a
principally polarised $g$-dimensional abelian variety as it has a fixed point on
Siegel's upper half space. Such an abelian variety $A$ together with an action
of $G$ may be seen as a principally polarised abelian variety over the
classifying stack $BG$ of $G$ and hence corresponds to a map $BG \to
\Acal_g$. Over the complex numbers, this map induces the given map $G \to
Sp_{2g}(\ZZ)$ on fundamental groups but the map makes sense over an arbitrary
base field given a principally polarised abelian $g$-dimensional variety with a
$G$-action and induces a pull back on cohomology.
%% \noindent
%% Let $G$ be a finite subgroup of ${\rm Sp}(2g,\ZZ)$. Every such
%% group is contained in the isotropy group of a point of the Siegel upper
%% half space
%% ${\Hcal}_g$ under the standard action of ${\rm Sp}(2g,\RR)$. Since
%% ${\Acal}_g$ is homotopy equivalent to
%% ${\rm Sp}(2g,\ZZ)$ we obtain a map of classifying spaces
%% $$
%% BG \longrightarrow {\Acal}_g.
%% $$
%% Over ${\Acal}_g$ we have the de Rham bundle $H_{dR}^1= R^1\pi_*\CC$,
%% where $\pi: {\Xcal}_g \to {\Acal}_g$ is the universal family (in the
%% orbifold sense). The de Rham bundle can be also obtained as the bundle
%% associated to the standard representation $\rho$  of ${\rm Sp}(2g,\ZZ)$:
%% it is the quotient of ${\rm Sp}(2g,\CC)\times V$ under
%% $(gp,\rho(p)^{-1}v)\sim (g,v)$ for all $g \in {\rm Sp}(2g,\CC)$ and $p \in
%% P$, where $P$ is a maximal parabolic subgroup (so that ${\Hcal}_g = {\rm
%% Sp}(2g,\CC)/P$). By universality the Chern classes $r_i$
%% of the de Rham bundle
%% map to the Chern
%% classes of the $2g$-dimensional representation of $G$. By a result of
%% Quillen the Chern classes of a faithful representation generate a subring
%% of the cohomology ring  over which the cohomology ring is finite as a
%% module. This already implies that not all Chern classes can vanish.
%% 
%% On the other hand by Hodge theory we have an exact sequence
%% $$
%% 0 \to \EE \longrightarrow H_{dR}^1 \longrightarrow \EE^{\vee} \to 0,
%% $$
%% which implies that all odd Chern classes $r_i$ of $H_{dR}^1$ vanish.
%% 
%% \smallskip
%;;%%%Lower bound for order Proposition lower bound
\begin{proposition}\label{lower bound}
Assume that $i \le g$.
The order of $r_{2i}$ is divisible by $n_i/2$ over $\CC$. In general, 
each prime $\ell$ different from the characteristic of the base field, $r_{2i}$
in $\ell$-adic cohomology has order divisible by the $\ell$-part of $n_i/2$.
\end{proposition}
%% {\sl Proof.}
\begin{proof}
The $\ell$-adic part implies the integral cohomology part so we may pick a prime
$\ell$ different from the characteristic of the base field and look at $r_{2i}$
in $\ell$-adic cohomology. What we want to show is that if $\ell$ is odd and
$\ell-1|2i$ and $k$ is the largest integer such that $\ell^{k-1}|2i$, then
$r_{2i}$ has order at least $\ell^k$ and similarly for $\ell=2$. We will do this
by defining a map $B\ZZ/\ell^k \to \Acal_g$, such that inverse image of the
$r_{2i}$ to $H^{4i}(\ZZ/\ell^k,\ZZ_\ell)$ has order $\ell^k$. Now,
$H^{4i}(\ZZ/\ell^k,\ZZ_\ell)$ is isomorphic to $\ZZ/\ell^k$ and hence an element
in it has order $\ell^k$ if and only if its reduction modulo $\ell$ is
non-zero. As $H^{4i}(\ZZ/\ell^k,\ZZ_\ell)/\ell$ injects into
$H^{4i}(\ZZ/\ell^k,\ZZ/\ell)$ it is enough to show that the pulback of $r_i$ is
non-zero in $H^{4i}(\ZZ/\ell,\ZZ/\ell)$. Assume first that $\ell$ is
odd. Consider now any Galois cover $C \to \PP^1$ with Galois group
$\ZZ/\ell^k$ which is ramified at $0$, $1$ and $\infty$ with ramification group
of order $\ell^k$, $\ell^k$, and $\ell$ respectively (the existence of such a
cover follows directly from Kummer theory). By the Hurwitz formula the genus of
such a covering fulfills the relation
$2g-2=-2\ell^k+2(\ell^k-1)+\ell^{k-1}(\ell-1)$, i.e.,
$2g=\ell^{k-1}(\ell-1)$. The action of $\ZZ/\ell^k$ on $H^1(C,\QQ_\ell)$ has to
contain a copy of the irreducible representation given by the action of
$\ell^k$'th roots of unity as it must be faithful and as $H^1(C,\QQ_\ell)$ has
dimension $\ell^{k-1}(\ell-1)$ it must be equal to it. Furthermore, this action
gives a map $B\ZZ/\ell^k \to \Mcal_g \to \Acal_g$. If now $x \in
H^2(\ZZ/\ell,\ZZ/\ell)$ is the natural generator (the Bockstein of the generator
of $H^1(\ZZ/\ell,\ZZ/\ell)$ corresponding to the identity map $\ZZ/\ell \to
\ZZ/\ell$) then the total Chern class of it is $\prod_{(i,\ell)=1}(1+ix)$ which
equals $(1+x^{\ell-1})^{\ell^{k-1}}=1+x^{\ell^{k-1}(\ell-1)}$. This gives the
non-triviality when $g=\ell^{k-1}(\ell-1)/2$ and when $g > \ell^k(\ell-1)/2$ we
simply add a principally polarised factor on which $\ZZ/\ell^k$ acts trivially.

When $\ell=2$ we may assume that $k > 2$ as the lower bound to be
 proven for $k\leq 2$
is implied by the one for $k=3$. We then make essentially the same construction,
a Galois cover with group $\ZZ/2^{k}$ of $\PP^1$ ramified at three points with
ramification groups of order $2^k$, $2^k$, and $2$ respectively of genus $g=2^{k-3}$. The rest of the
argument is identical to the odd case.
\end{proof}
\begin{remark}
i) It follows from (\ref{fund rel}) that the $r_i$ are torsion 
already in the Chow groups. Our result gives a lower bound for this order but we 
don't know if this bound is sharp.

ii) From the complex point of view our geometric construction can be seen simply
as constructing an element of order $\ell^k$ in
$\mathrm{Sp}_{\ell^{k-1}(\ell-1)}(\ZZ)$. This can be done directly, in the odd
case one may consider the ring of $\ell^k$'th roots of unity $R=\ZZ[\zeta]$ with
the obvious action of $\ZZ/\ell^k$ and the symplectic form
$\langle\alpha,\overline{\beta}\rangle :=
\Tr(\alpha\beta(\zeta-\zeta^{-1})^{-\ell^{k}+\ell^{k-1}+1})$. This is obviously
a symplectic invariant form and that it is
indeed an integer-valued perfect pairing follows from the fact that the
different of $R$ is the ideal generated by $(\zeta-\zeta^{-1})^{\ell^{k}-\ell^{k-1}-1}$.

iii) When $g=\ell^k(\ell-1)/2$ we actually get a lower bound for the top Chern class
of the de Rham cohomology of the universal curve over $\Mcal_g$. However, there
is no direct analogue of the trick of adding a factor with trivial action so
this does not give a lower bound for all $g \ge \ell^k(\ell-1)/2$.
\end{remark}
\begin{theorem}
We have that $r_{2i+1}=0$ for all $i$ and that the order of $r_{2i}$ in integral
($\ell$-adic) cohomology equals (resp.\ the $\ell$-part of) $n_i/2$ for
$i \le g$.
\end{theorem}
\begin{proof}
This follows immediately from Props.\ \ref{lower bound} and \ref{upper bound}.
\end{proof}
\begin{corollary} The order of $\lambda_g$ is at least $n_g/2$.
\end{corollary}
%% \noindent
%% {\sl Proof.} 
\begin{proof}
The top Chern class of $H_{dR}^1$ is $\lambda_g^2$.
\end{proof}
\begin{remark}
Our upper and lower bounds for $r_{2i}$ are off by a (multiplicative) factor of
$2$. Furthermore, when $g=1$ the lower bound is the correct order.
\end{remark}
\end{section}
\begin{section}{A complete representative cycle for $\lambda_g$ on
${\Acal}_g\otimes \FF_p$}
\bigskip
Consider the closed algebraic subset $V_0$ of ${\Acal}_g\otimes \FF_p$ of
all abelian varieties with $p$-rank zero. By Koblitz 
(see [K])
we know that this is a pure codimension $g$ cycle on ${\Acal}_g\otimes
\FF_p$. It is a complete cycle since abelian varieties of $p$-rank $0$ cannot
degenerate. Any complete subvariety of ${\Acal}_g$ has codimension at least $g$
in ${\Acal}_g$, see \cite{vdG1} and \cite{O}. Let $\tilde{\Acal}_g$
be a toroidal compactification.  One of us has proved (see \cite{vdG1})

\begin{theorem} The cycle class of $V_0$ in $CH^g_{\QQ}(\tilde{\Acal}_g)$ 
 is given by the formula $[V_0]=(p-1)(p^2-1)\cdots (p^g-1)\, \lambda_g$.
\end{theorem}
 
This shows that a multiple of $\lambda_g$ is an \emph{effective}  cycle
in characteristic $p>0$.  It seems unknown whether a non-zero multiple of the
class $\lambda_g$ can be represented by an effective cycle in characteristic
zero.
\end{section}

\begin{section}{A representative cycle for $\lambda_g$ on the Satake compactification}

We let ${\Acal}_g^*$ be the minimal (`Satake')
compactification as defined in \cite{F-C}. The Satake compactification
${\Acal}_g^*$ is a disjoint union
$$
{\Acal}_g^* = {\Acal}_g \sqcup {\Acal}_{g-1}\sqcup \ldots \sqcup {\Acal}_0.
$$
There is a natural morphism $q: \tilde{\Acal}_g \to {\Acal}_g^*$ for
every toroidal compactification $\tilde{\Acal}_g$ of ${\Acal}_g$.
In [3] two representative cycles for the class $\lambda_g \in
CH^g_{\QQ}({\Acal}_g^*)$ were given. We recall the result.

\begin{lemma}  The class $q_*(\lambda_g)$ in $CH_{\QQ}^g({\Acal}_g^*)$ is
represented by a multiple of the fundamental class of the boundary
$B_g^*={\Acal}_g^* - {\Acal}_g$. \end{lemma}
%% \noindent
%% {\sl Proof.}
\begin{proof}
The class $\lambda_g$ vanishes up to torsion on ${\Acal}_g$;
for dimension reasons $q_*\lambda_g$ is then represented by a multiple of
the fundamental class of $B_g^*$.
\end{proof}
\noindent
\begin{proposition} The cycle class $[B_g^*]$ of
the boundary is the same in the Chow group $CH_{\QQ}^g({\Acal}_g^*)$ as a
multiple of the class of the closure of the 
image of ${\Acal}^*_{g-1}$ in ${\Acal}^*_g$ under the
map $[X] \mapsto [X \times E]$, with $E$ a generic elliptic curve.
\end{proposition}
%% \smallskip
%% \noindent
%% {\sl Proof.}
\begin{proof}
Consider (for $g>2$) the space ${\Acal}_{g-1,1}$ of products of
a principally polarized abelian variety  of dimension $g-1$ and an
 elliptic curve. It is the image of ${\Acal}_{g-1}  \times
{\Acal}_1$ in ${\Acal}_g$ under a morphism to ${\Acal}_g$ which 
 can be extended to a morphism
${\Acal}_{g-1}^* \times {\Acal}_1^* \to {\Acal}_g^*$.
Since ${\Acal}_1$ is the affine
$j$-line we find a rational
equivalence between the cycle class of a generic fibre ${\Acal}^*_{g-1}
\times \{ j \}$  with $j$ a generic point on the $j$-line and a multiple of the fundamental class of the boundary
$B_g^*$.   
%%$\square$
\end{proof}
\bigskip

\bigskip

\section{A representative cycle for $\lambda_g$ on ${\Acal}_g^{\prime}$}

There are several compactifications of ${\Acal}_g\otimes k$.  We first choose a
suitable smooth toroidal compactification $\tilde{\Acal}_g$ as constructed in
\cite{F-C}. We have a natural map $q: \tilde{\Acal}_g \to {\Acal}_g^*$. We call
the inverse image under $q$ of $\sqcup_{j\leq t} {\Acal}_{g-j}$ the moduli space
$\tilde{\Acal}_g^{(t)}$ of \Definition{rank $\leq t$ degenerations.}  This space
parametrizes semi-stable abelian varieties whose torus-part has rank $\leq t$.
Furthermore, we let ${\Acal}_g^{\prime}=\tilde{\Acal}_g^{(1)}$ be the moduli
space of rank $1$-degenerations, i.e.\ the inverse image of ${\Acal}_g \sqcup
{\Acal}_{g-1} \subset {\Acal}_g^*$ under the natural map $ q: \tilde{\Acal}_g
\to {\Acal}_g^*$. Unlike the higher rank space ${\Acal}_g^{(t)}$ for $t\geq 2$
the space ${\Acal}_g'$ does {\sl not} depend on a choice $\tilde {\Acal}_g$ of
compactification of ${\Acal}_g$; it  is a \emph{canonical}
partial compactification
on ${\Acal}_g$. If we want a full  compactification then there is not really
a unique one, but we must make choices.  See \cite{M2}.

We start by  applying the Grothendieck-Riemann-Roch theorem 
to the  structure sheaf $O_{\tilde{\Xcal}_g}$  of a
 compactification $\tilde{\Xcal}_g$ 
as constructed in [2],  of  the universal
semi-abelian variety and the
morphism $\pi: \tilde{\Xcal}_g \to \tilde{\Acal}_g$.
This gives in the Chow groups with rational coefficients
$$
ch(\pi_!O_{\tilde{\Xcal}_g})=\pi_*(e^0 \, {\rm Td}^{\vee}(
\Omega^1_{\tilde{\Xcal}_g/\tilde{\Acal}_g})).
$$
The relative cotangent sheaf fits in an exact sequence
$$
0 \to \Omega^1_{\tilde{\Xcal}_g/\tilde{\Acal}_g} \to \EE
\to {\Fcal} \to 0.
$$
with $\EE$ the Hodge bundle and ${\Fcal}$ a sheaf with support
where $\pi$ is not smooth.
Note that by \cite{F-C} we have
$$
\pi^*(\EE)= \Omega^1_{\tilde{\Xcal}_g}(\log)/
\pi^*(\Omega^1_{\tilde{\Acal}_g}(\log)),
$$
where $\log$ refers to logarithmic poles along the divisors at 
infinity. We get
$$
ch(\pi_!(O_{\tilde{\Xcal}_g})=\pi_*(F) {\rm Td}^{\vee}(\EE)
$$
with $F := \mathrm{Td}^{\vee}({\Fcal})^{-1}$.
The derived sheaf $\pi_!(O_{\tilde{\Xcal}_g})$ equals $\wedge^* \EE=
\sum_{i=0}^g (-1)^i \wedge^i \EE$. By the Borel-Serre formula
we have $ch(\wedge^* \EE)= \lambda_g {\rm Td}(\EE)^{-1}$. Comparing
the terms of degree $\leq g$ yields  the result  of \cite{vdG1} :

\begin{proposition}
We have $\pi_*(F)=\lambda_g$. 
\end{proposition}

 From now on we work on the moduli space ${\Acal}_g^{\prime}=
\tilde{\Acal}_g^{(1)}$ of rank $\leq 1$ degenerations.  Let $D^0$ be the closed
subset corresponding to rank 1 degenerations. The divisor $D^0$ has a morphism
to $\phi:D^0 \to {\Acal}_{g-1}$ which exhibits $D^0$ as a quotient of the
universal abelian variety over ${\Acal}_{g-1}$.  The fibre over $x \in
{\Acal}_{g-1}$ is the dual ${\hat X}_{g-1}$ of the abelian variety $X$
corresponding to $x$. The `universal' semi-abelian variety $G$ over ${\hat
{\Xcal}} _{g-1}$ is the $\GG_m$-bundle obtained from the Poincar\'e bundle $P
\to {\Xcal}_{g-1} \times {\hat {\Xcal} }_{g-1}$ by deleting the zero-section. We
have the maps
$$
G=P -\{ (0)\} \to {\Xcal}_{g-1} \times_{{\Acal}_g} \hat{\Xcal}_{g-1} 
\to {\Acal}_{g-1}.
$$
We shall now work out an expression for $\pi_*(F)$. Let $B_g$ be the
cycle on ${\Acal}_g^{\prime}$ which is the locus of {\sl trivial}
semi-abelian extensions
$$
1\to \GG_m \to \tilde{X} \to X_{g-1}\to 0
$$
The cycle $B_g$ sits in $D^0$ as the zero section of $\phi: D^0 \to
{\Acal}_{g-1}$.

\begin{theorem} The Chow class  of $B_g$ and that of
${\Acal}_{g-1} \times \{j\}$ in $CH_{\QQ}^g({\Acal}_g^{\prime})$ are both
equal to  $(-1)^g{\lambda_g}/{\zeta (1-2g)}$ with $\zeta (s)$ denoting the
Riemann  zeta function. 
\end{theorem}
\smallskip
\noindent

The fibre over $x\in D^0$ is a compactification $\bar{G}$ of a
$\GG_m$-bundle $G$ over an abelian variety $X_{g-1}$ of dimension
$g-1$. The points where $\pi$ is not smooth are exactly the
points of $\bar{G}-G$. So globally the locus where $\pi$ is not
smooth is a cycle $\Delta$ obtained from glueing the $0$-section
and the $\infty$-section of the $\PP^1$-bundle associated to the
Poincar\'e bundle $P$ over ${\Xcal}_{g-1}\times_{{\Acal}_{g-1}}
\hat{\Xcal}_{g-1}$. 
 Note that $\Delta$ is of codimension $2$ and that an \'etale cover of
 $\Delta$ is given by ${\Xcal}_{g-1} \times {\hat {\Xcal}}_{g-1}$.                          

The normal bundle to $\Delta$ on ${\Xcal}_{g-1} \times
\hat{\Xcal}_{g-1}$
is then $N= P \oplus \tau^*(P^{-1})$ with $\tau(x, \hat x)= x+\hat x,
\hat x)$.
(We identify $\hat X$ with $X$ if needed.)
We write $\alpha_1= c_1(P)$ and $\alpha_2= c_1(\tau^*(P^{-1}))$.
On the space of rank $\leq 1$ degenerations (an \'etale cover of which
is  ${\Acal}_g^{(1)}= {\Acal}_g
 \cup {\hat{\Xcal}}_{g-1}$) we can write $\Delta= \alpha_1\alpha_2$.
 Let $i: \Delta \to {\tilde{\Xcal}}_g$ be the inclusion.

Then if we write
$$
{\rm Td}^{\vee}(L)=\frac{c_1(L)}{ (e^{c_1(L)}-1)}=
\sum_{k=0}^{\infty} \frac{b_k}{ k!} t^k
$$
we have (cf.\  Mumford [10], p.\ 303):
$$
\pi_*({\rm Td}^{\vee}(O_{\Delta}^{-1} -1) =
\pi_*(\sum_{k=1}^{\infty} \frac{(-1)^kb_{2k}}{(2k)!} 
i_*\left(\frac{ \alpha^{2k-1} + \alpha_2^{2k-1}}{ \alpha_1 + \alpha_2}\right).
\eqno(2)
$$
Observe that $P|X_{g-1}\times \hat x= L_{\hat x}$ and
$\tau^*(P^{-1})|X_{g-1}\times \hat x = t_{\hat x}^*(L_{\hat x}^{-1})$. This
implies by the Theorem of the Square that $P \otimes \tau^*(P^{-1})|X_{g-1}\times \hat x$ is trivial, i.e.
we have
$$
c_1(N)= c_1(P\otimes \tau^*(P^{-1}))= \alpha_1+\alpha_2=\hat{p}^*(\beta)
$$
with $\beta$ a codimension $1$ class on ${\hat{\Xcal}}_{g-1}$
and $\hat{p}$ the projection of ${\Xcal}_{g-1}\times_{{\Acal}_{g-1}}
\hat{\Xcal}_{g-1}$ on the second factor.
In order to
determine $\beta$ we restrict to the other fibres ($x \times \hat X_{g-1}$)
and find (writing $p=p_X$ and $\hat p = p_{\hat X}$):
$$
P = L_x = t_x^*T\otimes T^{-1} 
$$
and

\begin{align} \notag
\tau^*(P^{-1})& =\tau^*(m^*T^{-1}\otimes p^*T \otimes {\hat p}^*T)\\
\notag 
&= 2^*t_x^*T^{-1} \otimes t_x^*T \otimes T\\ \notag
&= L_x^{-1}\otimes (t_x^*T^{-1} \otimes (-1)^*t_x^*(T)^{-1})\\
\notag
\end{align}
\noindent
We find (assuming that the $\Theta$-divisor  is symmetric)
$$
\beta= -2T \qquad {\rm on}\quad \hat X.
$$
and
$$
N= P \oplus P^{-1}\otimes ({\hat p}^*(O(-2T))).\eqno(3)
$$
We can consider this as a global identity on $\Delta$. The
line bundle $T$ restricts in each fibre $\hat{X}$ to the
theta divisor.
Developing the terms in (2) we get expressions

\begin{align} \notag
\pi_*(i_*(\alpha_1+\alpha\notag_2)^r(\alpha_1\alpha_2)^s))&=
\pi_*(i_*({\hat p}^*(\beta^r)(\alpha_1\alpha_2)^s)\notag
\\ &= j_*(\beta^r{\pi^{\prime}}_*(\Delta^s)),\notag \\ \notag
\end{align}
where $\pi^{\prime}$ is the restriction to the boundary of $\pi$ and $j:D \to
{\Acal}_g^{(1)}$ is the inclusion of the boundary and where
we use $\pi i = \hat{p}$. 

For dimension reasons the only surviving terms are of the form
$j_*(\beta^r){\pi^{\prime}}_*(\Delta^{g-1})$. Thus the only term that
contributes is:
$$
\frac{(-1)^g b_{2g}}{ (2g)! }
\pi_*i_*((-1)^{g-1}(2g-1)(\alpha_1\alpha_2)^{g-1}).
$$
So we need to compute ${\pi^{\prime}}_*i_*(\Delta^{g-1})= \pi_*(\Delta^g)$.
The identity $\pi_*({\rm Td}^{\vee}({\Fcal})^{-1})=\lambda_g$ with
$$
F= {\rm Td}^{\vee}({\Fcal})^{-1}= 1+ f_2+f_4 +\ldots
$$
implies that
$$
\pi_*(f_{2g})= \frac{(-1) b_{2g} }{ (2g)! }
\pi_*i_*((2g-1)(\alpha_1\alpha_2)^{g-1})=\lambda_g.
$$
Represent $P$ by the divisor $\Pi$ on
${\Xcal}_{g-1} \times_{{\Acal}_{g-1}} \hat{\Xcal}_{g-1} $. Then we have by (3)
$$
\alpha_1^{g-1}\alpha_2^{g-1}
=(-1)^{g-1}( \Pi^{2g-2} + \sum_{j=1}^{g-2} {g-1 \choose j}
\Pi^{g-1+j}{\hat p}^*(2T)^{g-1-j} ).
$$
Now  apply  GRR to the  bundle $P\otimes {\hat p}^*(O(nT))$ on
${\Xcal}_{g-1} \times_{{\Acal}_{g-1}} \hat{\Xcal}_{g-1}$
; it says 
$$
ch({\hat p}_!(P\otimes {\hat p}^*O(nT)))= {\hat p}_*(e^{\Pi}) \, 
 \otimes e^{nT}.
$$ 
But ${\hat p}_!(P\otimes O(nT))$ is a sheaf with support (in  codimension
$g-1$) over the zero section $S$. Once again 
by applying GRR, this time to the inclusion $S \to \hat{\Xcal}_{g-1}$,
we get the relation
 $$c_{g-1}({\hat p}_!(P\otimes
O(nT))= (-1)^{g-2}(g-2)!\,  S
$$
It follows that 

$${\hat p}_*(\Pi^{2g-2-j})T^j  = 0\quad {\rm   if} \quad  j\neq 0$$
and
$${\hat p}_*(\Pi^{2g-2-j})T^j=(-1)^{g-1}(2g-2)!\, [S] \quad
{\rm if } \quad   j=0.$$
So we find
$$
\pi_*(\Delta^g)=j_*({\hat p}_*(\Pi^{2g-2}))= (-1)^{g-1} (2g-2)! B_g,
$$
where $B_g$ is the zero section of ${\Xcal}_{g-1} \to {\Acal}_{g-1}$.

We thus get

\begin{align} \notag
\lambda_g & = \frac{(-1) b_{2g} }{(2g)! }
\pi_*i_*((2g-1)(\alpha_1\alpha_2)^{g-1}) \\ \notag
& = \frac{(-1) b_{2g} } { (2g)! }(2g-1) (2g-2)! (-1)^{g-1} B_g \\ \notag
&= (-1)^g \zeta(1-2g)B_g\\ \notag
\end{align}

\begin{corollary} On the space of rank $\leq 1$ degenerations
${\Acal}_g^{\prime}$ we have the formula 
$\lambda_g= (-1)^g \zeta(1-2g)\, B_g $
with $B_g$ the locus of semi-stable abelian varieties which
are trivial extensions of an abelian variety with $\GG_m$.
\end{corollary}
\bigskip

\noindent
\begin{example}
We have $[B_1]= 12\, \lambda_1$ and $[B_2]= 120\, \lambda_2$.
\end{example}
\end{section}
\newcommand\eprint[1]{Eprint:~\texttt{#1}}
%\eject

%
%

\end{document}